# A blind spot in undergraduate mathematics: The circular definition of the length of the circle, and how it can be turned into an enlightening example

Alexei Vernitski, Department of Mathematical Sciences, University of Essex, Colchester, UK.

## Abstract

We highlight the fact that in undergraduate calculus, the number pi is defined via the length of the circle, the length of the circle is defined as a certain value of an inverse trigonometric function, and this value is defined via pi, thus forming a circular definition. We present a way in which this error can be rectified. We explain that this error is instructive and can be used as an enlightening topic for discussing different approaches to mathematics with undergraduate students.

## 1. Context

What I describe is borne out of many years of my experience of teaching first-year undergraduate mathematics. Also my approach benefits from my research interests. After a period of conducting research in pure mathematics, I diversified in two directions. On the one hand, I train artificial intelligence to produce proofs of simple mathematical results. On the other hand, I explore how producing and noticing shapes and patterns helps the brain to learn mathematics. In the former context, mathematics is "more formal" than usual, and in the latter context, mathematics is much "less formal" than usual. This dual perspective gives me a unique vantage point to appreciate the interplay between different approaches to mathematics. The purpose of this note is to highlight one example from first-year university mathematics which, as I will show, can become a topic for an eye-opening discussion on how "more formal" and "less formal" approaches to mathematics are both important for learning and understanding.

This duality of mathematics is perceived by many maths lecturers and, unfortunately, is felt as something which must be, at best, overcome or, at worst, reluctantly tolerated. A passionately written example of it is Allenby's book for first-year undergraduates *Numbers and Proofs*. Allenby writes about "good old days", when "one learnt to construct proofs", whereas "in recent years" "finding of strongly suggestive 'patterns' seems to have replaced *real* mathematical activity". Having said this, Allenby seems grudgingly to recognise persuasive limitations of "real" mathematics, saying that at the basic level, learners of mathematics should "be happy to accept certain easily believed, simple assertions as being unquestionably true" and at the advanced level, mathematicians "often have to read over proofs of results they 'know' (that is, 'feel in their bones') are true, just to see if the argument given is sufficient to establish a result they, in any case, believe".

However, as we discuss below, it is perfectly possible to appreciate the importance of both processes of "constructing proofs" and "feeling in one's bones that a result is true", and I believe it is possible and desirable to talk about this to students, using suitable examples like the one below.

## 2. The vicious circle

How do you define number $\pi$? If you are like me, you define it as a half of the length of the unit circle (and this definition, or an equivalent definition, is given in mathematical dictionaries). The next question is, how do we define the length of the unit circle? Undergraduate calculus (or engineering mathematics) tells us that the length of a curve represented by a function $y = f(x)$ is defined as $\int_a^b \sqrt{1 + \left(\frac{dy}{dx}\right)^2}\, dx$. Consider the unit circle defined by the equation $x^2 + y^2 = 1$. Let us use the definition of the length of the curve as expressed by the integral above to find the length of the top-right quarter of this circle. The function is $f(x) = \sqrt{1-x^2}$, and the limits or integration are $a = 0$, $b = 1$. Accordingly, the length of the quarter-circle is expressed as

$\int_0^1 \sqrt{1 + \frac{x^2}{1-x^2}}\, dx$, which can be simplified to $\int_0^1 \frac{dx}{\sqrt{1-x^2}}$. We can look up this integral in the table of indefinite integrals; it is an inverse trigonometric function $\sin^{-1}(x)$. To finish finding the value of the definite integral, we will need the value of $\sin^{-1}(x)$ at $x = 1$, so we look up the definition of $\sin^{-1}(x)$ and find out that the value of $\sin^{-1}(1)$ is defined to be equal to $\frac{\pi}{2}$. We have come full circle: unexpectedly, in undergraduate mathematics the value of number $\pi$ is defined via the value of number $\pi$.

One can attempt to bypass this circular definition by representing the circle in the parametric form. As we will see now, this does not help. Let us start from the beginning again; we want to define number $\pi$; it is defined as a half of the length of the unit circle; but what is the length of the unit circle? If a curve is given by a parametric representation $x = x(t)$, $y = y(t)$ then its length is defined as $\int_a^b \sqrt{\left(\frac{dx}{dt}\right)^2 + \left(\frac{dy}{dt}\right)^2}\, dt$. For the unit circle defined by the equation $x^2 + y^2 = 1$, according to undergraduate calculus (or engineering mathematics), we have $x = \cos t$, $y = \sin t$, and $t$ varies from $a = 0$ to $b = 2\pi$. Therefore, the integral turns into $\int_0^{2\pi} 1\, dt$ and is equal to $2\pi$. If we knew what the value of $\pi$ is, this would be a perfectly satisfactory answer. However, we do not know yet what $\pi$ is; we were taking this integral hoping to use the answer as a definition of $\pi$, and it does not help us, because it only says, basically, that the value of $2\pi$ is defined as being equal to $2\pi$.

Why is it a problem? This is a circular definition. Perhaps the value of $\pi$ actually does not exist. Perhaps function $\sin^{-1}(x)$ is not defined at $x = 1$. Perhaps the unit circle does not have length. By defining these concepts one via another we might be creating an illusion of a waterproof mathematical exposition, but what if none of these things exist?

The circular definitions in this section are not artificially constructed by me to bewilder you; as you can check by perusing textbooks on calculus or engineering mathematics, they really form a part of undergraduate courses in mathematics. As an aside, we can conjecture that if these circular arguments feature in published textbooks, they are even more likely to be present in individual lecturers' lecture notes.

Thus, we have to conclude that undergraduate mathematics leaves the question of the existence of $\pi$ open. Perhaps $\pi$, and the length of the circle, and the inverse trigonometric functions exist, or perhaps they do not.

## 3. How can an undergraduate student define $\pi$?

A good news is that $\pi$, and $\sin^{-1}(x)$, and the length of the circle exist. Here is how we could define them without circular definitions in a way which a determined undergraduate student could follow.

As we consider integral $\int_0^1 \frac{dx}{\sqrt{1-x^2}}$, of course we cannot say that it is $\sin^{-1}(x)$ because we don't know yet if function $\sin^{-1}(x)$ exists. However, with some effort we can write the Maclaurin series of a function whose derivative is $\frac{1}{\sqrt{1-x^2}}$; let us denote the series by $M(x)$. Unfortunately, it is difficult to prove that $M(x)$ converges, and it is especially difficult to prove that $M(1)$ converges, and I would not expect a first-year student to be much interested in all details of these proofs. But fortunately, and importantly, these proofs involve only real analysis and number theory and don't depend on the existence of $\pi$.

After we have satisfied ourselves that one can prove that $M(x)$ converges, recall that $\int_0^1 \frac{dx}{\sqrt{1-x^2}} = M(1)$ is the length of a quarter of the unit circle; therefore, the length of the unit circle exists and is equal to $4 \int_0^1 \frac{dx}{\sqrt{1-x^2}} = 4M(1)$. Number $\pi$ is a defined as a half of this length, that is, $\pi$ exists and is equal to $2M(1)$, or, to be more specific, $\pi = 2 + \frac{1}{3} + \frac{3}{20} + \cdots$.

At last, we can define $\sin^{-1}(x)$ to be equal to $M(x)$. Note that function $\sin x$ is still not defined. Indeed, $\sin x$ is defined in undergraduate calculus (or engineering mathematics) either using angles in triangles or using arc lengths in the unit circle. While we did not know what $\pi$ is, we could not measure angles in radians, and while we did not know that the unit circle has length, we could not measure arc lengths. Now that we know that function $\sin^{-1}(x)$ exists, we can define $\sin x$ as the inverse function of $\sin^{-1}(x)$.

There can be other ways in which we can remove circularity in the definitions described in the previous section, but the one I describe in the previous paragraphs seems to me the most practical one in the context of the undergraduate curriculum. One other way to define $\pi$ is as twice the smallest positive solution of the equation $cos\ x = 0$, where function $cos\ x$ is defined as a Maclaurin series. Remmert states that "the definition of $\frac{1}{2}\pi$ as the smallest positive zero of $cos\ x$ is now commonplace", although in my experience, I have never seen this definition (nor any other non-circular definition) in textbooks on calculus or engineering mathematics.

## 4. Turning it into an enriching experience

Now that I know that the definition of $\pi$ and the length of the circle in undergraduate mathematics is circular and, therefore, invalid, what am I supposed to do about it?

One approach is to ignore this fact and to continue teaching in the usual way, without telling my students about it. However, if I do this, I miss an opportunity to talk to my students about interesting and important details of the way mathematics is developed and used.

Another approach is to stress the validity of the intuitive understanding of mathematics. If we equip ourselves with a measuring tape and examples of circles (for examples, some pots and pans), we can convince ourselves that circles have length, and $\pi$ exists. Then I can tell students that the formal definition of the circle length and $\pi$ in undergraduate mathematics is flawed, but it does not matter, because from our experience we know that circles have length, and $\pi$ exists. Although this approach is practical, I believe it sends a wrong message to students; we want to feel that mathematics is robust and logical.

Yet another approach is to shun the intuitive understanding of $\pi$ as invalid, and follow the formal approach thoroughly, defining the value of $\pi$ to be $2 + \frac{1}{3} + \frac{3}{20} + \cdots$, as I described in the previous section. Although there is some appeal in this approach, I believe it is unnatural. For instance, we want to be able to use the usual geometric definition of $\pi$. The series $2 + \frac{1}{3} + \frac{3}{20} + \cdots$ looks contrived, and speaking practically, it is not even a good way to approximate the value of $\pi$. We want to be able to use a geometric definition of $\sin x$, rather than defining it as an inverse function of a certain Maclaurin series.

I believe that the best way is to combine approaches above and to explain to students different facets of mathematical practice. It is true that our experiments with pots and pans do enable us to produce perfectly workable definitions of the circle length and $\pi$. I can tell my students that it is an inspiring example of ingenuity of human mind that after some experiments with circles, in one stroke of genius we can conceive of a range of useful mathematical concepts. It is true that in comparison with this intuitive ingenuity, formal mathematics feels like a blunt tool; I can tell my students that a mathematician needs to write hundreds of pages defining limits, series, integrals etc. before she is able to define the circle length and $\pi$. However, on the positive side, we must recognise that the mathematician can achieve the result, and the circle length and $\pi$ can be defined formally. Therefore, separately from our intuitive ingenuity at conceiving mathematical concepts, we should also celebrate the mastery of mathematicians at wielding formal mathematical tools. As we do this, we should also appreciate the surprising reach of formal mathematics; formal mathematical arguments might feel circuitous, but they eventually catch up with our expectations, enriching our understanding on the way.

What I described is a rich picture of mathematics; however, this is not yet the whole story. I can also tell my students that formal mathematics is intentionally moulded to fit in with our informal understanding of mathematics. Indeed, mathematicians actively wanted circles to be defined in such a way that they have length, and $\pi$ to be defined so that its value matches the value measured in experiments. This is why mathematicians spent two thousand years perfecting definitions of limits, series, integrals etc. until these definitions, whereas satisfactorily formal and logical, also matched our informal concepts of the circle length and $\pi$.

Understanding and untangling the circular definition described in this note, and making it formal, involves the use of a wide range of topics of undergraduate mathematics. This, in itself, is enough to make this circular definition worth discussing with students. However, what makes this circular definition a special example to me is how disengaged our intuitive understanding of the circle length and $\pi$ and the formal definitions of them are from each other. This gap is the reason why this circular definition appears repeatedly in calculus (or engineering mathematics) textbooks and remains unnoticed by students. This is where an opportunity arises; this circular definition is a perfect material on which the role of different approaches to mathematics can be explored and discussed.